\documentclass{article}
\usepackage[utf8]{inputenc}
\usepackage{amsthm}
\usepackage{amsmath, amssymb, amsthm, mathtools}
\usepackage{graphicx}
\usepackage{bbm}
\usepackage{comment}
\usepackage[margin=1.5in]{geometry}
\usepackage{tikz}
\usepackage{mathtools}
\usepackage{authblk}
\usepackage{hyperref}

\DeclareMathOperator{\N}{\mathbf{N}}

\DeclareMathOperator{\R}{\mathbf{R}}
\DeclareMathOperator{\E}{\mathbb{E}}
\DeclareMathOperator{\Z}{\mathbf{Z}}
\DeclareMathOperator{\V}{\mathbb{V}}
\DeclareMathOperator{\PP}{\mathbb{P}}

\def\longto{\longrightarrow}
\def\disp{\displaystyle}

\newtheorem{theorem}{Theorem}[section]

\newtheorem{proposition}[theorem]{Proposition}

\theoremstyle{remark}

\theoremstyle{remark}

\newenvironment{customproof}[1][Proof]{%
  \begin{proof}[\textnormal{\textbf{#1}}]%
}{%
  \end{proof}%
}

\title{Asymptotic Results for Uniform Group Drawing in the Coupon Collector's Problem}
\author[1,2]{Daniel Berend}
\author[2]{Tomer Sher}
\affil[1]{Institute for the Theory of Computing, Ben-Gurion University, Beer Sheva 84105, Israel.}
\affil[2]{Department of Mathematics, Ben-Gurion University, Beer Sheva 84105, Israel.}

\graphicspath{ {./images/} }

\begin{document}

\maketitle

\begin{abstract} 

The article explores the asymptotic behavior of the expected number of drawings in the Coupon Collector’s Problem with group-drawing under the uniform distribution. In this variant, each draw consists of a package of $s$ distinct coupons selected uniformly at random from a set of $n$ coupons. We focus on three regimes of the package size $s$: (i) constant $s$, (ii) $s$ proportional to $n$, and (iii) $s$ ``very close'' to $n$. For each case, we provide precise asymptotic expressions for the expected collection time. 

\end{abstract}

\section{Introduction}


Consider the problem of collecting a complete set of $n$ distinct coupons, where each cereal box contains one coupon chosen uniformly at random from the $n$ types. The Coupon Collector Problem (CCP) asks: what is the expected number of boxes one may expect to purchase to collect all $n$ coupons? This fundamental problem arises in probability theory, computer science, and combinatorial analysis, with applications ranging from algorithm design to data structure analysis. Notably, the problem can be traced back to de Moivre~\cite{tod}.

It is well known, and easy to prove, that the expected number of drawings needed to complete the collection is
\begin{align} \label{1}
    n\cdot H_n = n\cdot\left(1+\frac{1}{2}+\cdots+\frac{1}{n-1}+\frac{1}{n}\right),
\end{align}
where $H_n$ is the $n$-th harmonic number.

Let $W_{n,M}$ be the number of drawings until an $M$-collection (that is, $M$ copies of each coupon) is obtained. Erd\H{o}s and R\'{e}nyi~\cite{erdos} found the asymptotic distribution of $W_{n,M}$ as $n\to\infty$ for fixed $M$:
\begin{align*}
    \lim_{n \to \infty} \PP\left(W_{n,M} \le n\left(\log n+(M-1)\log \log n-\log (M-1)!+x\right)\right) 
    = \PP(X \le x), \qquad x \in \mathbb{R},
\end{align*}
for $X \sim \text{Gumbel}(0,1)$.
Recall~\cite{Gumbel1954,PinheiroFerrari} that $X \sim \text{Gumbel}(\mu,\beta)$ with parameters $\mu \in \mathbb{R}$ and $\beta > 0$ if its distribution function is
\[
F(x) = \exp\left(-e^{-(x-\mu)/\beta}\right), \qquad x \in \mathbb{R}.
\]

Coupon subset collection is a natural generalization of the classical version. Instead of selecting a single coupon each time, one selects a set of $s$ distinct coupons (henceforth a \emph{package}), where $s$ is a fixed integer between $1$ and $n-1$. We are interested in the distribution of the number of rounds needed to obtain all coupons.  
This batched-draw variant also arises in biological systems. For example, Giannakis et~al.~\cite{giannakis2022exchange} model mitochondrial DNA exchange in plant cells as a network-based coupon collection process, where encounters correspond to group draws of genetic elements.

Schilling and Henze~\cite{henze} extended the Erd\H{o}s--R\'enyi asymptotic result to this setting. Let $W_{n,s,M}$ denote the number of drawings needed to obtain an $M$-collection when each package consists of a set of $s$ coupons, where $s$ is a constant and all $\binom{n}{s}$ such sets are equally likely. They proved that
\begin{align*}
    \lim_{n \to \infty} 
    \PP\left( W_{n,s,M} \le \frac{n}{s}
    \left(\log n+(M-1)\log \log n-\log (M-1)!+x\right)\right)
    = \PP(X \le x), 
    \qquad x \in \mathbb{R},
\end{align*}
for $X \sim \text{Gumbel}(0,1)$.

\medskip

While the asymptotic behavior of $W_{n,s,M}$ is known for constant $s$, less is known when $s$ varies with $n$. 
In particular, the expected collection time may exhibit different qualitative behaviors depending on the relationship between $s$ and $n$. In this paper, we study the asymptotics of this expectation for three typical regimes of $s$ -- small, moderate, and large (see next section for the exact regimes).

In Section~\ref{main results} we state our main results.  
The proofs are presented in Section~\ref{proofs}.

 \section{The Main Results}\label{main results}

Let $Y_{n,s}$ the number of draws required to collect all $n$ coupons when the distribution over all $s$-subsets is uniform. Usually, we will simply write $Y$, leaving $n$ and $s$ implicit.

We study the asymptotic behavior of $Y$ in some special cases:
\begin{itemize}
    \item{} Case I. $s$ is a constant.
    \item{} Case II. $s=c\cdot n$ for an arbitrary fixed $0<c<1$.
    \item{} Case III. $s = n-\lambda\cdot n^\beta$ for an arbitrary fixed $0< \beta <1$ and $\lambda >0$.
\end{itemize}
To be precise, in Cases II and III we need to round the right-hand side, but we will ignore it.
The three following theorems relate to these three cases, respectively.
\begin{theorem} \label{expectation s constant}
     If $s \ge 2$ is a constant, then
     \begin{align*}
         \E[Y] = \frac{n}{s}\cdot H_{n} -\frac{s-1}{2s}\cdot H_{n} + O\left(1\right),
     \end{align*}
      as $n\to\infty$.
\end{theorem}

At first glance, one may think that, as the expected value for $s=1$ is $n H_n$, for general constant $s$ it should be $n H_n /s$. The reason for its being slightly less is that collecting groups gives a slight advantage. The first coupon in the group may be any of the $n$ coupons, but the second is guaranteed to be different from the first, the third is guaranteed to be different from the first two, and so forth. This reduces the expected collection time a little for constant $s$, and by more and more as $s$ grows.

\begin{figure}[!ht]
\centering
\includegraphics[width=1\textwidth]{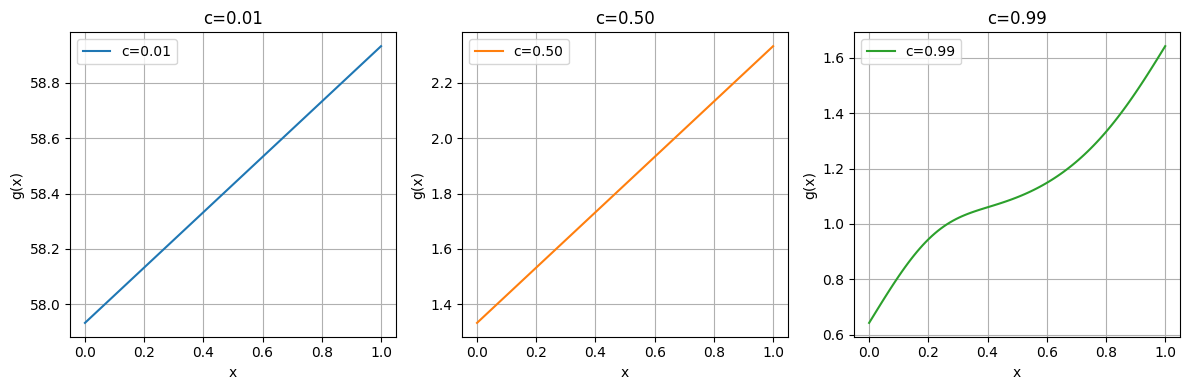}
\caption{$g(x)$ for $c = 0.01, 0.5, 0.99$}
\label{g(x)}
\end{figure}

We pass to Case II. Let $\{y\}$ denote the fractional part of a real number $y$. Define a function $g:[0,1)\longto\R$ by:
$$g(x)=\sum_{i=1}^{\infty} \left(1-e^{-(1-c)^{i-x-1}}-e^{-(1-c)^{-i-x}}\right), \qquad x\in [0,1).$$
Note that $g$ is continuous and bounded in $[0,1)$. In Figure \ref{g(x)} we plot $g$ on $[0,1)$ for three values of $c$: $0.01, 0.5, 0.99$. Note that only for $c$ close to~1 can we distinguish $g$ from a linear function with slope 1.

\begin{theorem} \label{expectation s=cn}
    Let $s = c\cdot n$, where $0 < c < 1$ is a constant, and let $\alpha=\alpha(n) = \{\log_{1/(1-c)} n\}$. As $n\to\infty$, we have
    \begin{equation}\label{caseII_asymptotics}
    \E[Y] = \lfloor  \log_{1/(1-c)} n \rfloor + g(\alpha(n)) + o(1).
    \end{equation}  
\end{theorem}
The first term on the right-hand side of (\ref{caseII_asymptotics}) is the main term. The other two jointly may be replaced by $O(1)$, but give more information as they are written. Notice that $g$ grows by~1 as we move from $0$ to $1^-$.
Thus, the right-hand side of (\ref{caseII_asymptotics}), considered as a function of a real variable $x$ instead of $n$, changes continuously as $x$ grows from $1/(1-c)^{m-\varepsilon}$ to $1/(1-c)^{m+\varepsilon}$ for a positive integer $m$.

In Case III we will distinguish between the subcase where $\beta$ is a rational of the form $\beta = (t-1)/t$ for some integer $t\ge 2$, and the subcase where it is not. Theorem \ref{expectation fractional} describes the asymptotic number of drawings in each subcase.

A sequence $(X_{n})_{n=1}^{\infty}$ of random variables \emph{converges in distribution} to a random variable~$X$, and we write $X_{n} \xrightarrow[n\to \infty]{\mathcal{D}} X$, if $\lim_{n \to \infty} F_{n}(x) = F(x)$ for all $x \in \mathbb{R}$ at which $F$ is continuous, where $F$ and $F_{n}$ denote the  cumulative distribution functions of $X$ and $X_{n}$, respectively. It \emph{converges in probability}, and we write $X_{n} \xrightarrow[n\to \infty]{\mathcal{P}} X$, if $\lim_{n\to\infty}\mathbb{P}\big(|X_n-X| > \varepsilon\big) = 0$ for all $\varepsilon > 0$.
Recall that, for two sequences $a_n$ and $b_n$ such that $b_n>0$ for sufficiently large $n$, we write $a_n = \omega(b_n)$ if $\lim_{n \to \infty}\frac{a_n}{b_n} = \infty$.

\begin{theorem} \label{expectation fractional}
     Let $t\ge 2$ be an integer. \newline\newline
    (a) If $s = n -\lambda\cdot n^{\frac{t-1}{t}}$, where $\lambda >0$, then $$\lim_{n \to \infty} \E[Y] = te^{-\lambda^t}+(t+1)(1-e^{-\lambda^t}).$$ Moreover, $Y \xrightarrow[n\to \infty]{\mathcal{D}} Z$, where $Z$ is a random variable with the distribution $$\PP[Z=i] = 
    \begin{cases}
            e^{-\lambda^t}, &i=t, \\
            1-e^{-\lambda^t}, &i=t+1.
        \end{cases}$$
    (b) If $s = n-o(n^{t/(t+1)})$ and $s = n-\omega(n^{(t-1)/t})$ then 
    $$ \lim_{n \to \infty} \E[Y] = t+1.$$
    Moreover, $Y \xrightarrow[n\to \infty]{\mathcal{P}} t+1$.
\end{theorem}

\section{Proofs}\label{proofs}

First we deal with the proof of case I, i.e., where $s$ is a constant.

\begin{customproof}[Proof of Theorem \ref{expectation s constant}]

Since $\frac{1}{1-x} = 1 + x + O(x^2)$ for $|x|<1$, for a fixed positive $a$:
\begin{align*}
    \frac{1}{n-a} = \frac{1}{n} \cdot \frac{1}{1-\frac{a}{n}} = \frac{1}{n}\cdot\left(1 + \frac{a}{n} + O\left(\frac{1}{n^2} \right) \right).
\end{align*}
Hence,
\begin{align*}
    H_n - H_{n-s} &= \frac{1}{n} +  \frac{1}{n-1} + \cdots + \frac{1}{n-s+1} \\
    &= \frac{1}{n}\cdot\left[\left(1+\frac{0}{n}\right) + \left(1+\frac{1}{n}\right)+ \cdots + \left(1+\frac{s-2}{n}\right) + \left(1+\frac{s-1}{n}\right) + s\cdot O\left(\frac{1}{n^2} \right)  \right] \\
    &= \frac{1}{n}\cdot\left[s + \frac{s(s-1)/2}{n} + O\left(\frac{1}{n^2} \right)  \right] \\&= \frac{s}{n}\cdot\left[1 + \frac{s-1}{2n} + O\left(\frac{1}{n^2} \right) \right].
\end{align*}
It follows that
\begin{align*}
    \frac{H_n}{H_n - H_{n-s}} &=  \frac{\frac{n}{s}H_n}{1 + \frac{s-1}{2n} + O\left(\frac{1}{n^2} \right)} \\&= \frac{n}{s}H_n\cdot\left(1 - \frac{s-1}{2n} + O\left(\frac{1}{n^2} \right) \right) \\&= \frac{n}{s}H_n - \frac{s-1}{2s}\cdot H_n + o(1).
\end{align*}
By \cite[Prop. 2.3]{Berend&Sher}, we obtain
$$\E[Y] = \frac{n}{s}H_n - \frac{s-1}{2s}\cdot H_n + O(1).$$ 
\end{customproof}

For cases II and III, we will employ the following result \cite[Thm. 8.3.1]{alon}.
\newline\newline
\textbf{Theorem A.}
    \textit{Let X = X(n) be a sum of indicator random variables, and let $\mu >0$. If $\lim_{n \to \infty} \E[\binom{X}{r}] = \mu^r/r!$ for every $r$, then $X \xrightarrow[n\to \infty]{\mathcal{D}}$ \rm{Po}$(\mu)$.}
\newline

It will be convenient to deal first with case III. We prove Theorem \ref{expectation fractional} by using Proposition \ref{distribution of X_t} below. We have $n$ coupons and the package size is $s = n -\lambda\cdot n^{\frac{t-1}{t}}$, where $t \ge 2$ is an integer and $\lambda > 0$. Denote by $X_{k}$ the number of those coupons not obtained during the first $k$ steps of the process, $k=1,2,3,\ldots$.
Notice that $X_{k} = \sum_{i=1}^n \mathbbm{1}^k_{i}$, where $\mathbbm{1}^k_{i}$ is the indicator of the event whereby the $i$-th coupon has not been collected within the first $k$ steps. 
\begin{proposition} \label{distribution of X_t}
    Let $s = n -\lambda\cdot n^{\frac{t-1}{t}}$, where $t \ge 2$ is an arbitrary fixed integer and $\lambda > 0$. Then
    \begin{align*}
        X_{t} \xrightarrow[n\to \infty]{\mathcal{D}} \mathrm{Po}(\lambda^{t}).
    \end{align*}
\end{proposition}

\begin{customproof} 
For any fixed $r$, the variable $\binom{X_{k}}{r}$ is a sum of indicators: $$\binom{X_{k}}{r} = \sum_{1\le i_{1}<\cdots<i_{r}\le n} \mathbbm{1}^k_{i_{1},\ldots,i_{r}},$$
where $\mathbbm{1}^k_{i_{1},\ldots,i_{r}}$ is the indicator of the event whereby all the $r$ coupons $i_{1},\ldots,i_{r}$ have not been collected within the first $k$ steps. Denote by $A_{i}$ the event whereby the $i$-th coupon has not been collected within the first $t$ steps. Then:
\begin{align*}
    \E \binom{X_{t}}{r} &= \sum_{1\le i_{1}<\cdots<i_{r}\le n} \E[\mathbbm{1}^t_{i_{1},\ldots,i_{r}}] = \binom{n}{r}\cdot \PP(A_{1} \cap A_{2} \cap \cdots \cap A_{r}) \\&= \frac{n!}{r!\cdot(n-r)!}\cdot \left(\frac{(n-r)(n-r-1)\cdot\cdots\cdot(n-r-s+1)}{n(n-1)\cdot\cdots\cdot (n-s+1)}\right)^t \\&= \frac{1}{r!}\cdot\frac{[(n-s)\cdot(n-s-1)\cdot\cdots\cdot(n-s-r+1)]^t}{[n\cdot(n-1)\cdot\cdots\cdot(n-r+1)]^{t-1}} \\&= \frac{1}{r!}\cdot\frac{[\lambda n^{\frac{t-1}{t}}\cdot(\lambda n^{\frac{t-1}{t}}-1)\cdot\cdots\cdot(\lambda n^{\frac{t-1}{t}}-r+1)]^t}{[n\cdot(n-1)\cdot\cdots\cdot(n-r+1)]^{t-1}} \xrightarrow[n \to \infty]{} \frac{\lambda^{tr}}{r!}.
\end{align*}

By Theorem A, we have $X_{t} \xrightarrow[n\to \infty]{\mathcal{D}} \mathrm{Po}(\lambda^{t})$.
\end{customproof}

With this proposition in hand, we now prove Theorem \ref{expectation fractional}.

\begin{customproof}[Proof of Theorem \ref{expectation fractional}]
(a) First we prove the ``moreover'' part. We will show that $\lim_{n \to \infty} \PP[Y = i] =0$ for $0 \le i \le t-1$ and for $i\ge t+2$. We bound the probability $\PP[Y \ge i]$. Denote by $B_{j}$ the event whereby the $j$-th coupon has not been collected within $i-1$ steps. Then:
\begin{align} \label{prob upper bound}
    \PP[Y \ge i] &= \PP\left[\bigcup_{j=1}^n B_{j}\right] \le \sum_{j=1}^n \PP[B_{j}] = n\cdot\left(\frac{n-s}{n}\right)^{i-1} = \lambda^{i-1}\cdot n^{1-\frac{i-1}{t}}.
\end{align}
Hence $\lim_{n \to \infty} \PP[Y = i] =0$  for $i \ge t+2$. 

We will show that the same holds for $0 \le i \le t-1$ as well. To this end, we first calculate the expectation and variance of $X_{t-1}$. Using the representation of $X_{t-1}$ as a sum of indicators, we obtain

\begin{align*}
    \E[X_{t-1}] &= \sum_{i=1}^n \E[\mathbbm{1}^{t-1}_{i}] = n\cdot\left(\frac{n-s}{n}\right)^{t-1} = \frac{\lambda^{t-1}n^{(t-1)^2/t}}{n^{t-2}} =\lambda^{t-1}n^{1/t},
\end{align*}
and
\begin{align*}
    \E[X_{t-1}^2] &= \sum_{i=1}^n \E[\mathbbm{1}^{t-1}_{i}] + \sum_{1 \le i \neq j \le n} \E[\mathbbm{1}^{t-1}_{i}\cdot \mathbbm{1}^{t-1}_{j}] \\&=  n\cdot\left(\frac{n-s}{n}\right)^{t-1} + (n^2-n)\cdot\left(\frac{(n-s)(n-s-1)}{n(n-1)}\right)^{t-1} \\&= \lambda^{t-1}n^{1/t} + \lambda^{2t-2}n^{1/t}\cdot(n-1)\cdot\left(\frac{n^{(t-1)/t}-1}{n-1}\right)^{t-1}.
\end{align*}
Thus,
\begin{align*}
    \V[X_{t-1}] &= \lambda^{t-1}n^{1/t} + \lambda^{2t-2}n^{1/t}\cdot(n-1)\cdot\left(\frac{n^{(t-1)/t}-1}{n-1}\right)^{t-1} - \lambda^{2t-2}n^{2/t} \\&= \lambda^{t-1}n^{1/t} + \lambda^{2t-2}n^{2/t}\cdot\frac{(1-n^{1/t}/n)^{t-1}}{(1-1/n)^{t-2}} - \lambda^{2t-2}n^{2/t}\\&= \lambda^{t-1}n^{1/t} + \lambda^{2t-2}n^{2/t}\cdot(1+o(1)) - \lambda^{2t-2}n^{2/t}\\&= \lambda^{t-1}n^{1/t} + \lambda^{2t-2}n^{2/t}\cdot o(1) = o(n^{2/t}).
\end{align*}
By the second moment method,
\begin{align*}
    \PP[Y \le t-1] &= \PP[X_{t-1} = 0] \\&\le \frac{\V[X_{t-1}]}{\E[X_{t-1}]^2} \\&= \frac{\V[X_{t-1}]}{\lambda^{2t-2}n^{2/t}} \xrightarrow[n \to \infty]{} 0.
\end{align*}
We conclude that $ \lim_{n \to \infty} \PP[Y = i] =0$ for $0 \le i \le t-1$. By Proposition \ref{distribution of X_t}, 
\begin{align*}
    \PP[Y = t] = \PP[Y \le t]-\PP[Y \le t-1] = \PP[X_{t} = 0] -o(1) \xrightarrow[n \to \infty]{} e^{-\lambda^t}.
\end{align*}
Therefore:
\begin{align*}
    \lim_{n \to \infty} \PP[Y=t+1]  = 1-e^{-\lambda^t}. 
\end{align*}

To evaluate $\E[Y]$ asymptotically, we will obtain lower and upper bounds. The lower bound is trivial: 
\begin{align} \label{exp lower bound}
  \E[Y]\ge t\cdot\PP[Y=t]+(t+1)\cdot\PP[Y=t+1]. 
\end{align}
For the upper bound we use the Tail-Sum Formula \cite[pp.215--220]{tail-sum} and (\ref{prob upper bound}),

\begin{align*} 
    \E[Y] &= \sum_{i=1}^\infty \PP[Y \ge i] = \sum_{i=1}^{t+1} \PP[Y \ge i] + \sum_{i=t+2}^\infty \PP[Y \ge i] \\&\leq \sum_{i=1}^{t+1} \PP[Y \ge i] + \sum_{i=t+2}^\infty n^{(t+1)/t}\cdot\left(\frac{\lambda}{n^{1/t}}\right)^i \\&= t\cdot\PP[Y=t]+(t+1)\cdot\PP[Y=t+1] +o(1) + n^{(t+1)/t}\cdot\sum_{i=t+2}^\infty \left(\frac{\lambda}{n^{1/t}}\right)^i \\&= t\cdot\PP[Y=t]+(t+1)\cdot\PP[Y=t+1] + o(1) + n^{(t+1)/t}\cdot\frac{\lambda^{t+2}\cdot n^{-(t+2)/t}}{1-\lambda\cdot n^{-1/t}} \\&= t\cdot\PP[Y=t]+(t+1)\cdot\PP[Y=t+1] +o(1), 
\end{align*} 
which yields:
\begin{align} \label{exp upper bound}
    \E[Y] \le t\cdot\PP[Y=t]+(t+1)\cdot\PP[Y=t+1] +o(1). 
\end{align}
By (\ref{exp lower bound}) and (\ref{exp upper bound}):

\begin{align*}
    \lim_{n \to \infty} \E[Y] = te^{-\lambda^t}+(t+1)(1-e^{-\lambda^t}).
\end{align*}
(b) Since $s = n-o(n^{t/(t+1)})$ and $s = n-\omega(n^{(t-1)/t})$, for every $\lambda_{1},\lambda_{2} > 0$ we have $$ n-\lambda_{2}\cdot n^{(t-1)/t} < s < n-\lambda_{1}\cdot n^{t/(t+1)} $$ for sufficiently large $n$. For simplicity, denote $s'=n-\lambda_{2}\cdot n^{(t-1)/t}$ and $s''=n-\lambda_{1}\cdot n^{t/(t+1)}$. Obviously, $Y_{s''} \preceq Y_{s} \preceq Y_{s'}$, where $X_1 \preceq X_2$ means that $X_1$ is stochastically smaller than~$X_2$. In particular, $\E[Y_{s''}] \le \E[Y_{s}] \le \E[Y_{s'}]$ as well. We have an upper bound for $\E[Y_{s'}]$ and a lower bound for $\E[Y_{s''}]$ from the proof of part (a). This yields

\begin{align*}
    (t+1)\PP[Y_{s''} = t+1]+(t+2)\PP[Y_{s''} = t+2] \le \E[Y_{s}]
\end{align*}
and
\begin{align*}
    \E[Y_{s}] \le  t\PP[Y_{s'} = t]+(t+1)\PP[Y_{s'} = t+1] + o(1).
\end{align*}
Letting $n \to \infty$, we obtain 
\begin{align*}
    (t+1)e^{-\lambda_{1}^{t+1}}+(t+2)(1-e^{-\lambda_{1}^{t+1}}) \le 
 	\varliminf_{n \to \infty} \E[Y_{s}]
\end{align*}
and
\begin{align*} 
    \varlimsup_{n \to \infty} \E[Y_{s}] \le te^{-\lambda_{2}^t}+(t+1)(1-e^{-\lambda_{2}^t}).
\end{align*}
Taking $\lambda_{1} \to 0$ and $\lambda_{2} \to \infty$, we obtain
\begin{align*}
    \lim_{n \to \infty} \E[Y_{s}] = t+1.
\end{align*}
Now we prove that $Y_{s} \xrightarrow[n\to \infty]{\mathcal{P}} t+1$. Since $Y_{s} \preceq Y_{s'}$, we have $\PP[Y_{s'} \le t+1] \le \PP[Y_{s} \le t+1]$. By (a), $$1-o(1) \le \PP[Y_{s} \le t+1] \le 1,$$ and therefore $$\lim_{n \to \infty} \PP[Y_{s} \le t+1] = 1.$$ On the other hand, $Y_{s''} \preceq Y_{s}$, so that $\PP[Y_{s} \le t] \le \PP[Y_{s''} \le t]$. By (a), $$ 0 \le \PP[Y_{s} \le t] \le o(1),$$ and therefore $$\lim_{n \to \infty} \PP[Y_{s} \le t] = 0.$$ In conclusion,
\begin{align*}
    \lim_{n \to \infty} \PP[Y_{s} = t+1] = 1.
\end{align*}
\end{customproof}

Finally, we deal with case II. For the proof of Theorem \ref{expectation s=cn}, we start with

\begin{proposition} \label{distribution of X_k+i}
    Let $s=c\cdot n$ for some $0<c<1$. Denote by $X_{\ell}$ the number of unseen coupons after $\ell$ steps. Let $k(n)=  \lfloor  \log_{1/(1-c)} n \rfloor$ and $\alpha(n) = \{\log_{1/(1-c)} n\}$. For any increasing sequence $(n_{j})_{j=1}^\infty$ of positive integers, if $$\lim_{j \to \infty} \alpha( n_{j}) =  \alpha, \quad (\alpha \in [0,1]),$$ then $$X_{k(n_{j}) + i} \xrightarrow[j\to \infty]{\mathcal{D}} \text{\rm{Po}}((1-c)^{i-\alpha}), \quad i \in \Z.$$
\end{proposition}

\begin{customproof}
First, we calculate $\E \binom{X_{k(n)+i}}{r}$ in order to use Theorem A. Similarly to the proof of Proposition \ref{distribution of X_t}, denote by $A_{l}$ the event whereby the $l$-th coupon has not been collected within the first $k(n)+i$ steps. Writing $\binom{X_{k(n)+i}}{r}$ as a sum of indicators, we obtain:
\begin{align*}
    \E \binom{X_{k(n)+i}}{r} &= \binom{n}{r}\cdot \PP(A_{1} \cap A_{2} \cap \cdots \cap A_{r}) \\&= \frac{n!}{r!\cdot(n-r)!}\cdot \left(\frac{(n-r)(n-r-1)\cdot\cdots\cdot(n-r-s+1)}{n(n-1)\cdot\cdots\cdot(n-s+1)}\right)^{k(n)+i} \\&= \frac{1}{r!}\cdot\frac{[(1-c)n\cdot((1-c)n-1)\cdot\cdots\cdot((1-c)n-r+1)]^{k(n)+i}}{[n\cdot(n-1)\cdot\cdots\cdot(n-r+1)]^{k(n)+i-1}} \\&= \frac{1}{r!}\cdot n(n-1)\cdot\cdots\cdot(n-r+1)\\&\quad\cdot\left(\frac{(1-c)n\cdot((1-c)n-1)\cdot\cdots\cdot((1-c)n-r+1)}{n\cdot(n-1)\cdot\cdots\cdot(n-r+1)}\right)^{k(n)+i} \\&= \frac{1}{r!}\left[n\cdot\left(\frac{(1-c)n}{n}\right)^{k(n)+i}\right]\cdot\left[(n-1)\cdot\left(\frac{(1-c)n-1}{n-1}\right)^{k(n)+i}\right]\\&\quad\cdot\cdots\cdot\left[(n-r+1)\cdot\left(\frac{(1-c)n-r+1}{n-r+1}\right)^{k(n)+i}\right].
\end{align*}
For simplicity, we denote $d = 1/(1-c)$:
\begin{align*}
    \E \binom{X_{k(n)+i}}{r} &= \frac{1}{r!}\left[n\cdot\left(\frac{(1-c)n}{n}\right)^{\log_{d}n+i-\alpha(n)}\right]\cdot\left[(n-1)\cdot\left(\frac{(1-c)n-1}{n-1}\right)^{\log_{d}n+i-\alpha(n)}\right]\\&\quad\cdot\cdots\cdot\left[(n-r+1)\cdot\left(\frac{(1-c)n-r+1}{n-r+1}\right)^{\log_{d}n+i-\alpha(n)}\right] \\&= \frac{1}{r!}\left[n\cdot\left(\frac{(1-c)n}{n}\right)^{\log_{d}n+i-\alpha(n)}\right]\cdot\left[(n-1)\cdot\left(\frac{(1-c)n-1}{n-1}\right)^{\log_{d}n+i-\alpha(n)}\right]\\&\quad\cdot\cdots\cdot\left[(n-r+1)\cdot\left(\frac{(1-c)n-r+1}{n-r+1}\right)^{\log_{d}n+i-\alpha(n)}\right] \\&=\frac{1}{r!}\left[n\cdot\left(\frac{(1-c)n}{n}\right)^{\log_{d}n}\right]\cdot\left[(n-1)\cdot\left(\frac{(1-c)n-1}{n-1}\right)^{\log_{d}n}\right]\\&\quad\cdot\cdots\cdot\left[(n-r+1)\cdot\left(\frac{(1-c)n-r+1}{n-r+1}\right)^{\log_{d}n}\right]\\&\quad\cdot\left(\frac{1}{d^r}\cdot\frac{(n-d)\cdot\cdots\cdot(n-d(r-1))}{(n-1)\cdot\cdots\cdot(n-r+1)}\right)^{i-\alpha(n)}.
\end{align*}
By the previous calculation,
 
\begin{align*}
    \E \binom{X_{k(n_{j})+i}}{r} &= \frac{1}{r!}\underbrace{\left[n_{j}\cdot\left(\frac{(1-c)n_{j}}{n_{j}}\right)^{\log_{d}n_{j}}\right]}_{\xrightarrow[j \to \infty]{} 1}\cdot\underbrace{\left[(n_{j}-1)\cdot\left(\frac{(1-c)n_{j}-1}{n_{j}-1}\right)^{\log_{d}n_{j}}\right]}_{\xrightarrow[j \to \infty]{} 1}\cdot\\&\quad\cdot \cdots\cdot\underbrace{\left[(n_{j}-r+1)\cdot\left(\frac{(1-c)n_{j}-r+1}{n_{j}-r+1}\right)^{\log_{d}n_{j}}\right]}_{\xrightarrow[j \to \infty]{}1}\cdot\\&\quad\cdot\underbrace{\left(\frac{1}{d^r}\cdot\frac{(n_{j}-d)\cdot\cdots\cdot(n_{j}-d(r-1)}{(n_{j}-1)\cdot\cdots\cdot(n_{j}-r+1)}\right)^{i-\alpha(n_{j})}}_{\xrightarrow[j \to \infty]{} (1-c)^{r\cdot(i-\alpha)}}.
\end{align*}
Therefore,
\begin{align*}
    \lim_{j \to \infty} \E \binom{X_{k(n_{j})+i}}{r} = \frac{(1-c)^{(i-\alpha)r}}{r!},
\end{align*}
and by Theorem A, $$X_{k(n_{j}) + i} \xrightarrow[j \to \infty]{\mathcal{D}} \text{Po}\left((1-c)^{i-\alpha}\right).$$ 
\end{customproof}

    As in (\ref{prob upper bound}),
    \begin{align} \label{prob upper bound for s=cn}
        \PP[Y \ge i] \le n\cdot\left(\frac{n-s}{n}\right)^{i-1} = (1-c)^{i-1-\log_{1/(1-c)}n}= (1-c)^{i-1-k(n)-\alpha(n)}.
    \end{align}
    Now let $l \in \N$. Using again the Tail-Sum Formula, (\ref{prob upper bound for s=cn}), and Proposition \ref{distribution of X_k+i}, we obtain the upper bound
    \begin{align*}
        \E[Y] &= \sum_{i=1}^\infty \PP[Y \ge i] \\&=  \sum_{i=1}^{k(n_{j})-l} \PP[Y \ge i] + \sum_{i=k(n_{j})-l+1}^{k(n_{j})+l} \PP[Y \ge i] + \sum_{i=k(n_{j})+l+1}^\infty \PP[Y \ge i] \\&\le  (k(n_{j})-l) +\sum_{i=-l+1}^{l} \PP[Y \ge k(n_{j})+i] + \sum_{i=k(n_{j})+l+1}^\infty (1-c)^{(i-1)-k(n_{j})-\alpha(n_{j})} \\&= (k(n_{j})-l) +\sum_{i=-l+1}^{l} \PP[X_{k(n_{j})+i-1} > 0] + \sum_{i=0}^\infty(1-c)^{i+l-\alpha(n_{j})} \\&= (k(n_{j})-l) +\sum_{i=-l+1}^{l} \left(1-e^{-(1-c)^{i-\alpha-1}}+o(1)\right) + \frac{(1-c)^{l-\alpha(n_{j})}}{c} \\&= (k(n_{j})-l) +\sum_{i=1}^{l} \left(1-e^{-(1-c)^{i-\alpha-1}}\right) +\sum_{i=0}^{-(l-1)} \left(1-e^{-(1-c)^{i-\alpha-1}}\right) \\&\quad +\frac{(1-c)^{l-\alpha(n_{j})}}{c} +o_{l}(1) \\&= (k(n_{j})-l) +\sum_{i=1}^{l} \left(1-e^{-(1-c)^{i-\alpha-1}}\right) +\sum_{i=1}^{l} \left(1-e^{-(1-c)^{-i-\alpha}}\right) \\&\quad +\frac{(1-c)^{l-\alpha(n_{j})}}{c} +o_{l}(1) \\&= k(n_{j}) -\sum_{i=1}^{l} e^{-(1-c)^{-i-\alpha}} + \sum_{i=1}^{l} \left(1-e^{-(1-c)^{i-\alpha-1}}\right) +\frac{(1-c)^{l-\alpha(n_{j})}}{c} + o_{l}(1).
    \end{align*}
    Letting $j \to \infty$, we obtain
    \begin{align*}
        \varlimsup_{j \to \infty} \E[Y-k(n_{j})] \le \sum_{i=1}^{l} \left(1-e^{-(1-c)^{i-\alpha-1}}\right) -\sum_{i=1}^{l} e^{-(1-c)^{-i-\alpha}} +\frac{(1-c)^{l-\alpha}}{c}.
    \end{align*}
    Letting now $l \to \infty$:
    \begin{align*}
         \varlimsup_{j \to \infty} \E[Y-k(n_{j})] \le \sum_{i=1}^{\infty} \left(1-e^{-(1-c)^{i-\alpha-1}} - e^{-(1-c)^{-i-\alpha}}\right).
    \end{align*}
    To get a lower bound, denote by $X^i_{m}$ the number of times the $m$-th coupon has been seen within $i$ steps. 
    We have:
    \begin{align*}
        \PP[Y \le i] &= \PP\left[X^{i}_{1} \ge 1, X^{i}_{2} \ge 1,\ldots,X^{i}_{n} \ge 1\right] \\&= \PP\left[X^{i}_{1} \ge 1\right] \cdot \PP\left[X^{i}_{2} \ge 1 | X^{i}_{1} \ge 1\right] \cdot\cdots\cdot \PP\left[X^{i}_{n} \ge 1 | X^{i}_{1} \ge 1, \ldots,  X^{i}_{n-1} \ge 1 \right].
    \end{align*}
    Consider a typical factor $\PP\left[X^{i}_{m} \ge 1 | X^{i}_{1} \ge 1, \ldots,  X^{i}_{m-1} \ge 1\right]$ on the right-hand side. The conditioning event $\{X^{i}_{1} \ge 1, \ldots,  X^{i}_{m-1} \ge 1\}$ gives positive information about appearances of some coupons other than coupon $m$, and therefore the conditional probability is at most $\PP\left[X^{i}_{m-1} \ge 1\right]$. It follows that: 
    \begin{align*}
        \PP[Y \le i] &\le \prod_{j=1}^n \PP[X^{i}_{j} \ge 1] = \prod_{j=1}^n  \left(1-\left(\frac{n-s}{n}\right)^{i}\right) = (1-(1-c)^i)^n.
    \end{align*}
    \newline
    For simplicity, in the next calculation we denote $x= e^{-1/(1-c)^{k(n_j)+\alpha(n_j)+1}}$. Let $m_{0}$ be the
    smallest integer such that $1/(1-c)^{m_{0}} > 2$. Similarly to the upper bound, using the Tail-Sum Formula, and Proposition \ref{distribution of X_k+i}, we have the lower bound

    \begin{align*}
        \E[Y] &= \sum_{i=1}^\infty \PP[Y \ge i] \\&=  \sum_{i=1}^{k(n_{j})-l} \PP[Y \ge i] + \sum_{i=k(n_{j})-l+1}^{k(n_{j})+l} \PP[Y \ge i] + \sum_{i=k(n_{j})+l+1}^\infty \PP[Y \ge i] \\&\ge \sum_{i=1}^{k(n_{j})-l} (1-\PP[Y < i])+\sum_{i=-l+1}^{l} \PP[Y \ge k(n_{j})+i] \\&= (k(n_{j})-l) -\sum_{i=0}^{k(n_{j})-l-1} \PP[Y \le i]+\sum_{i=-l+1}^{l} \PP[X_{k(n_{j})+i-1} > 0] \\&= (k(n_{j})-l)+\sum_{i=-l+1}^{l} (1-e^{-(1-c)^{i-\alpha-1}}+o(1)) - \sum_{i=0}^{k(n_{j})-l-1} \PP[Y \le i] \\&\ge (k(n_{j})-l)+ \sum_{i=1}^{l} \left(1-e^{-(1-c)^{i-\alpha-1}}\right) +\sum_{i=1}^{l} \left(1-e^{-(1-c)^{-i-\alpha}}\right) \\&\quad -\sum_{i=0}^{k(n_{j})-l-1} (1-(1-c)^i)^{n_{j}} +o_{l}(1).
    \end{align*}
    
    Therefore,
    \begin{align*}
         \E[Y] &\ge (k(n_{j})-l) +l+\sum_{i=1}^{l} (1-e^{-(1-c)^{i-\alpha-1}}) - \sum_{i=1}^{l} e^{-(1-c)^{-i-\alpha}} \\&\quad -\sum_{i=0}^{k(n_{j})-l-1}e^{-n_{j}(1-c)^i} +o_{l}(1) \\&=k(n_{j}) +\sum_{i=1}^{l} (1-e^{-(1-c)^{i-\alpha-1}}) - \sum_{i=1}^{l} e^{-(1-c)^{-i-\alpha}}-\sum_{i=1}^{k(n_{j})-l} e^{-n_{j}(1-c)^{i-1}} +o_{l}(1)\\&= k(n_{j})-\sum_{i=1}^{l} e^{-(1-c)^{-i-\alpha}} + \sum_{i=1}^{l} \left(1-e^{-(1-c)^{i-\alpha-1}}\right)-\sum_{i=1}^{k(n_{j})-l} x^{(1-c)^i} +o_{l}(1)\\&\ge k(n_j) -\sum_{i=1}^{l} e^{-(1-c)^{-i-\alpha}} + \sum_{i=1}^{l} \left(1-e^{-(1-c)^{i-\alpha-1}}\right) \\&\quad -2x^{(1-c)}\left(1-1/2^{k(n_j)-l}\right) +o_{l}(1).
    \end{align*}
    Letting $j \to \infty$, we obtain
    \begin{align*}
        \varliminf_{j \to \infty} \E[Y-k(n_j)] \ge \sum_{i=1}^{l} \left(1-e^{-(1-c)^{i-\alpha-1}}\right)-\sum_{i=1}^{l} e^{-(1-c)^{-i-\alpha}}. 
    \end{align*}
    Now let $l \to \infty$:
    \begin{align*}
         \varliminf_{j \to \infty} \E[Y-k(n_j)] \ge \sum_{i=1}^{\infty} \left(1-e^{-(1-c)^{i-\alpha-1}}-e^{-(1-c)^{-i-\alpha}}\right).
    \end{align*}
    Finally,
    \begin{align*}
        \lim_{n \to \infty} \E[Y-k(n_j)] = \sum_{i=1}^{\infty} \left(1-e^{-(1-c)^{i-\alpha-1}}-e^{-(1-c)^{-i-\alpha}}\right).
    \end{align*}
    \end{customproof}
    
\end{comment}

\begin{customproof}[Proof of Theorem \ref{expectation s=cn}]

The proof relies on a detailed asymptotic analysis of the tail-sum representation of the expectation. 
The process is decomposed into three parts -- contributions from drawings far below, near, and far above the critical index.
Using precise probabilistic estimates (mainly Proposition~\ref{distribution of X_k+i}), the main contribution is shown to come from a bounded window around the critical index, while the other parts vanish asymptotically. 
Upper and lower bounds are derived by applying exponential approximations and a second-moment argument, and the two bounds coincide in the limit, yielding the desired asymptotic formula.

Denote $k(n) = \lfloor  \log_{1/(1-c)} n \rfloor$ as in Proposition \ref{distribution of X_k+i}. Let $\ell$ be an arbitrary fixed positive integer. Using the Tail-Sum Formula, we decompose $\E[Y]$ as follows:
\begin{equation}\label{split}
  \begin{array}{rcl}
  \E[Y] & = & \disp\sum_{i=1}^\infty \PP[Y \ge i] \\\\
    & = & \disp\sum_{i=1}^{k(n_{j})-\ell} \PP[Y \ge i] + \sum_{i=k(n_{j})-\ell+1}^{k(n_{j})+\ell} \PP[Y \ge i] + \sum_{i=k(n_{j})+\ell+1}^\infty \PP[Y \ge i].
    \end{array}
\end{equation}
We first bound $\E[Y]$ from above. The first term on the right-hand side of (\ref{split}) is bound trivially:
\begin{equation}\label{first-term}
  \disp\sum_{i=1}^{k(n_{j})-\ell} \PP[Y \ge i] \le k(n_j) - \ell.
\end{equation}
For the second term, by Proposition \ref{distribution of X_k+i}, we have:
\begin{equation}\label{second-term}
  \begin{array}{rcl}
     \disp\sum_{i=k(n_{j})-\ell+1}^{k(n_{j})+\ell} \PP[Y \ge i]  
     & = & \disp\sum_{i=-\ell+1}^{\ell} \PP[Y \ge k(n_{j})+i]  \\\\
     & = & \disp\sum_{i=-\ell+1}^{\ell} \PP[X_{k(n_{j})+i-1} > 0]  \\\\
     & = & \disp\sum_{i=-\ell+1}^{\ell} \left(1-e^{-(1-c)^{i-\alpha-1}}+o(1)\right)  \\\\
     & = & \disp\sum_{i=1}^{\ell} \left(1-e^{-(1-c)^{-i-\alpha}}\right) + \sum_{i=1}^{\ell} \left(1-e^{-(1-c)^{i-\alpha-1}}\right) + o_{\ell}(1) \\\\
     & = & \ell +\disp\sum_{i=1}^{\ell} \left(1-e^{-(1-c)^{i-\alpha-1}} - e^{-(1-c)^{-i-\alpha}}\right) + o_{\ell}(1).
\end{array}
\end{equation}
Regarding the third term, as in (\ref{prob upper bound}),
\begin{align} \label{prob upper bound for s=cn}
    \PP[Y \ge i] \le n\cdot\left(\frac{n-s}{n}\right)^{i-1} = (1-c)^{i-1-\log_{1/(1-c)}n}= (1-c)^{i-1-k(n)-\alpha(n)},
\end{align}
and therefore
\begin{equation}\label{third-term}
  \begin{array}{rcl}
     \disp\sum_{i=k(n_{j})+\ell+1}^\infty \PP[Y \ge i]  
     & \le & \disp\sum_{i=k(n_{j})+\ell+1}^\infty(1-c)^{i-1-k(n_j)-\alpha(n_j)}  \\\\
     & = & \disp\frac{(1-c)^{\ell-\alpha(n_{j})}}{c}.
\end{array}
\end{equation}
Plugging (\ref{first-term}), (\ref{second-term}), and (\ref{third-term}) in (\ref{split}), we obtain
\begin{align*}
    \E[Y] - k(n_{j}) \le \disp\sum_{i=1}^{\ell} \left(1-e^{-(1-c)^{i-\alpha-1}} - e^{-(1-c)^{-i-\alpha}} \right) +\frac{(1-c)^{\ell-\alpha(n_{j})}}{c} + o_{\ell}(1).
\end{align*}
Passing to the limit as $j \to \infty$, we obtain
\begin{align*}
    \varlimsup_{j \to \infty} (\E[Y]-k(n_{j})) \le \sum_{i=1}^{\ell} \left(1-e^{-(1-c)^{i-\alpha-1}} - e^{-(1-c)^{-i-\alpha}}\right) +\frac{(1-c)^{\ell-\alpha}}{c}.
\end{align*}
Now let $\ell \to \infty$:
\begin{equation}\label{less-equal}
     \varlimsup_{j \to \infty} (\E[Y]-k(n_{j})) \le \sum_{i=1}^{\infty} \left(1-e^{-(1-c)^{i-\alpha-1}} - e^{-(1-c)^{-i-\alpha}}\right).
\end{equation}
In the other direction, we start again at (\ref{split}). For the second term on the right we still use (\ref{second-term}). The third term is estimated trivially:
\begin{equation}\label{third-term-opposite}
     \disp\sum_{i=k(n_{j})+\ell+1}^\infty \PP[Y \ge i]  
      \ge 0.
\end{equation}
We turn to the first term on the right. To estimate $\PP[Y \ge i]$, decompose $X_i$ again in the form
$$X_i=\disp\sum_{t=1}^n X_{i,t},$$
where $X_{i,t}=1$ if coupon $t$ has not arrived within the first $i$ drawings and $X_{i,t}=0$ otherwise, $1\le t\le n$. The decomposition gives
$$\E[X_i]=(1-c)^i n.$$
Also, since the $X_{i,t}, 1\le t\le n,$ are negatively correlated, we have
$$\V[X_i]\le n \V[X_{i,1}] = (1-c)^i (1-(1-c)^i) n \le (1-c)^i n,$$
so that
$$\E[X_i^2] \le (1-c)^i n + (1-c)^{2i} n^2.$$
Employing the second-moment method (cf. \cite[Ch. 21]{jukna2011extremal}, we arrive at
$$\begin{array}{rcl}
P[X_i>0] & \ge & \dfrac{E^2[X_i])} {E[X_i^2]}\\\\
& \ge & \dfrac{(1-c)^{2i} n^2}{(1-c)^i n + (1-c)^{2i} n^2}\\\\
& = & \dfrac{(1-c)^{i} n}{1 + (1-c)^{i} n}\\\\
& \ge & 1- \dfrac{1}{(1-c)^{i} n}.
\end{array}$$
It follows that
\begin{equation}\label{first-term-opposite}
\begin{array}{rcl}
  \disp\sum_{i=1}^{k(n_{j})-\ell} \PP[Y \ge i] 
  & = & \disp\sum_{i=1}^{k(n_{j})-\ell} \PP[X_i >0]\\\\
  & \ge & \disp\sum_{i=1}^{k(n_{j})-\ell} \left(1- \dfrac{1}{(1-c)^{i} n}\right)\\\\
  & = & k(n_j)-\ell -\disp\sum_{i=1}^{k(n_{j})-\ell} \left(\dfrac{(1-c)^{k(n_j)+\alpha(n_j)}}{(1-c)^{i}}\right)\\\\
  & = & k(n_j)-\ell -\disp\sum_{i=1}^{k(n_{j})-\ell} (1-c)^{k(n_j)+\alpha(n_j) - i}\\\\
  & \ge & k(n_j)-\ell -\dfrac{(1-c)^{\ell+\alpha(n_j)}}{c}.
  \end{array}
\end{equation}
Plugging (\ref{first-term-opposite}), (\ref{second-term}), and (\ref{third-term-opposite}) in (\ref{split}), we obtain
\begin{align*}
    \E[Y] - k(n_{j}) \ge \disp\sum_{i=1}^{\ell} 
    \left(1-e^{-(1-c)^{i-\alpha-1}} - e^{-(1-c)^{-i-\alpha}} \right) +\frac{(1-c)^{\ell+\alpha(n_j)}}{c} + o_{\ell}(1).
\end{align*}
Taking limits as $j \to \infty$, we obtain
\begin{align*}
    \varliminf_{j \to \infty} (\E[Y]-k(n_{j})) \ge \sum_{i=1}^{\ell} \left(1-e^{-(1-c)^{i-\alpha-1}} - e^{-(1-c)^{-i-\alpha}}\right) -\dfrac{(1-c)^{\ell+\alpha}}{c}.
\end{align*}
Letting $\ell \to \infty$:
\begin{equation}\label{greater-equal}
     \varliminf_{j \to \infty} (\E[Y]-k(n_{j})) \ge \sum_{i=1}^{\infty} \left(1-e^{-(1-c)^{i-\alpha-1}} - e^{-(1-c)^{-i-\alpha}}\right).
\end{equation}
The two inequalities (\ref{less-equal}) and (\ref{greater-equal}) prove the theorem.
\end{customproof}

\section{Acknowledgments}
We thank the anonymous reviewers for their valuable comments and suggestions that helped improve this work.

\clearpage
\bibliographystyle{plain}
\bibliography{references_2}
\end{document}